\documentclass{amsart}

\input{epsf.tex}
\usepackage{amsfonts,amssymb,verbatim,amsmath,amsthm,latexsym,textcomp,amscd, hyperref}
\usepackage{latexsym,amsfonts,amssymb,epsfig,verbatim}
\usepackage{amsmath,amsthm,amssymb,latexsym,graphics,textcomp}
\usepackage{graphicx}
\usepackage{color}
\usepackage{url}

\begin{document}

\newtheorem{theorem}{Theorem}[section]
\newtheorem{prop}[theorem]{Proposition}
\newtheorem{lemma}[theorem]{Lemma}
\newtheorem{cor}[theorem]{Corollary}
\newtheorem{defn}[theorem]{Definition}
\newtheorem{conj}[theorem]{Conjecture}
\newtheorem{claim}[theorem]{Claim}
\newtheorem{rem}[theorem]{Remark}
\newtheorem{rmk}[theorem]{Remark}

\newcommand{\map}{\rightarrow}
\newcommand{\boundary}{\partial}
\newcommand{\C}{{\mathbb C}}
\newcommand{\integers}{{\mathbb Z}}
\newcommand{\natls}{{\mathbb N}}
\newcommand{\ratls}{{\mathbb Q}}
\newcommand{\reals}{{\mathbb R}}
\newcommand{\proj}{{\mathbb P}}
\newcommand{\lhp}{{\mathbb L}}
\newcommand{\tube}{{\mathbb T}}
\newcommand{\cusp}{{\mathbb P}}
\newcommand\AAA{{\mathcal A}}
\newcommand\BB{{\mathcal B}}
\newcommand\CC{{\mathcal C}}
\newcommand\DD{{\mathcal D}}
\newcommand\EE{{\mathcal E}}
\newcommand\FF{{\mathcal F}}
\newcommand\GG{{\mathcal G}}
\newcommand\HH{{\mathcal H}}
\newcommand\II{{\mathcal I}}
\newcommand\JJ{{\mathcal J}}
\newcommand\KK{{\mathcal K}}
\newcommand\LL{{\mathcal L}}
\newcommand\MM{{\mathcal M}}
\newcommand\NN{{\mathcal N}}
\newcommand\OO{{\mathcal O}}
\newcommand\PP{{\mathcal P}}
\newcommand\QQ{{\mathcal Q}}
\newcommand\RR{{\mathcal R}}
\newcommand\SSS{{\mathcal S}}
\newcommand\TT{{\mathcal T}}
\newcommand\UU{{\mathcal U}}
\newcommand\VV{{\mathcal V}}
\newcommand\WW{{\mathcal W}}
\newcommand\XX{{\mathcal X}}
\newcommand\YY{{\mathcal Y}}
\newcommand\ZZ{{\mathcal Z}}
\newcommand\CH{{\CC\HH}}
\newcommand\TC{{\TT\CC}}
\newcommand\EXH{{ \EE (X, \HH )}}
\newcommand\GXH{{ \GG (X, \HH )}}
\newcommand\GYH{{ \GG (Y, \HH )}}
\newcommand\PEX{{\PP\EE  (X, \HH , \GG , \LL )}}
\newcommand\MF{{\MM\FF}}
\newcommand\PMF{{\PP\kern-2pt\MM\FF}}
\newcommand\ML{{\MM\LL}}
\newcommand\PML{{\PP\kern-2pt\MM\LL}}
\newcommand\GL{{\GG\LL}}
\newcommand\Pol{{\mathcal P}}
\newcommand\half{{\textstyle{\frac12}}}
\newcommand\Half{{\frac12}}
\newcommand\Mod{\operatorname{Mod}}
\newcommand\Area{\operatorname{Area}}
\newcommand\ep{\epsilon}
\newcommand\hhat{\widehat}
\newcommand\Proj{{\mathbf P}}
\newcommand\U{{\mathbf U}}
 \newcommand\Hyp{{\mathbf H}}
\newcommand\D{{\mathbf D}}
\newcommand\Z{{\mathbb Z}}
\newcommand\R{{\mathbb R}}
\newcommand\Q{{\mathbb Q}}
\newcommand\E{{\mathbb E}}
\newcommand\til{\widetilde}
\newcommand\length{\operatorname{length}}
\newcommand\tr{\operatorname{tr}}
\newcommand\gesim{\succ}
\newcommand\lesim{\prec}
\newcommand\simle{\lesim}
\newcommand\simge{\gesim}
\newcommand{\simmult}{\asymp}
\newcommand{\simadd}{\mathrel{\overset{\text{\tiny $+$}}{\sim}}}
\newcommand{\ssm}{\setminus}
\newcommand{\diam}{\operatorname{diam}}
\newcommand{\pair}[1]{\langle #1\rangle}
\newcommand{\T}{{\mathbf T}}
\newcommand{\inj}{\operatorname{inj}}
\newcommand{\pleat}{\operatorname{\mathbf{pleat}}}
\newcommand{\short}{\operatorname{\mathbf{short}}}
\newcommand{\vertices}{\operatorname{vert}}
\newcommand{\collar}{\operatorname{\mathbf{collar}}}
\newcommand{\bcollar}{\operatorname{\overline{\mathbf{collar}}}}
\newcommand{\I}{{\mathbf I}}
\newcommand{\tprec}{\prec_t}
\newcommand{\fprec}{\prec_f}
\newcommand{\bprec}{\prec_b}
\newcommand{\pprec}{\prec_p}
\newcommand{\ppreceq}{\preceq_p}
\newcommand{\sprec}{\prec_s}
\newcommand{\cpreceq}{\preceq_c}
\newcommand{\cprec}{\prec_c}
\newcommand{\topprec}{\prec_{\rm top}}
\newcommand{\Topprec}{\prec_{\rm TOP}}
\newcommand{\fsub}{\mathrel{\scriptstyle\searrow}}
\newcommand{\bsub}{\mathrel{\scriptstyle\swarrow}}
\newcommand{\fsubd}{\mathrel{{\scriptstyle\searrow}\kern-1ex^d\kern0.5ex}}
\newcommand{\bsubd}{\mathrel{{\scriptstyle\swarrow}\kern-1.6ex^d\kern0.8ex}}
\newcommand{\fsubeq}{\mathrel{\raise-.7ex\hbox{$\overset{\searrow}{=}$}}}
\newcommand{\bsubeq}{\mathrel{\raise-.7ex\hbox{$\overset{\swarrow}{=}$}}}
\newcommand{\tw}{\operatorname{tw}}
\newcommand{\base}{\operatorname{base}}
\newcommand{\trans}{\operatorname{trans}}
\newcommand{\rest}{|_}
\newcommand{\bbar}{\overline}
\newcommand{\UML}{\operatorname{\UU\MM\LL}}
\newcommand{\EL}{\mathcal{EL}}
\newcommand{\tsum}{\sideset{}{'}\sum}
\newcommand{\tsh}[1]{\left\{\kern-.9ex\left\{#1\right\}\kern-.9ex\right\}}
\newcommand{\Tsh}[2]{\tsh{#2}_{#1}}
\newcommand{\qeq}{\mathrel{\approx}}
\newcommand{\Qeq}[1]{\mathrel{\approx_{#1}}}
\newcommand{\qle}{\lesssim}
\newcommand{\Qle}[1]{\mathrel{\lesssim_{#1}}}
\newcommand{\simp}{\operatorname{simp}}
\newcommand{\vsucc}{\operatorname{succ}}
\newcommand{\vpred}{\operatorname{pred}}
\newcommand\fhalf[1]{\overrightarrow {#1}}
\newcommand\bhalf[1]{\overleftarrow {#1}}
\newcommand\sleft{_{\text{left}}}
\newcommand\sright{_{\text{right}}}
\newcommand\sbtop{_{\text{top}}}
\newcommand\sbot{_{\text{bot}}}
\newcommand\sll{_{\mathbf l}}
\newcommand\srr{_{\mathbf r}}
\newcommand\geod{\operatorname{\mathbf g}}
\newcommand\mtorus[1]{\boundary U(#1)}
\newcommand\A{\mathbf A}
\newcommand\Aleft[1]{\A\sleft(#1)}
\newcommand\Aright[1]{\A\sright(#1)}
\newcommand\Atop[1]{\A\sbtop(#1)}
\newcommand\Abot[1]{\A\sbot(#1)}
\newcommand\boundvert{{\boundary_{||}}}
\newcommand\storus[1]{U(#1)}
\newcommand\Momega{\omega_M}
\newcommand\nomega{\omega_\nu}
\newcommand\twist{\operatorname{tw}}
\newcommand\modl{M_\nu}
\newcommand\MT{{\mathbb T}}
\newcommand\Teich{{\mathcal T}}
\renewcommand{\Re}{\operatorname{Re}}
\renewcommand{\Im}{\operatorname{Im}}

\title{Semiconjugacies Between Relatively Hyperbolic Boundaries}

\author{Shubhabrata Das}
\author{Mahan Mj}

\address{School of Mathematical Sciences, RKM Vivekananda University, Belur Math, WB-711 202, India}

\thanks{Research of first author partially supported by a Department of Science and Technology Research grant.
The second author is partly supported by a CSIR Junior Research Fellowship. This paper is part
of SD's PhD thesis  written under the supervision of MM. } 
\date{}

\begin{abstract}
We prove the existence of Cannon-Thurston maps for Kleinian groups corresponding to
 pared manifolds whose  boundary is incompressible away from cusps. We also describe the structure of these maps in terms of ending laminations.
\end{abstract}

\maketitle

\tableofcontents
\section{Introduction}
The aim of this paper is threefold: \\
1) To extend the main Theorems of \cite{mahan-split}, \cite{mahan-elct}   (which prove the existence  and structure
of Cannon-Thurston maps for surface groups without accidental parabolics) to  Kleinian groups corresponding to
{\it pared manifolds whose  boundary is incompressible away from cusps}. \footnote{A considerably more elaborate and somewhat
clumsier proof had been sketched in an earlier version of \cite{mahan-split}. This proof has been excised from the present version of \cite{mahan-split}.}
This is the content of Theorem \ref{main2}.\\
2) To give a considerably shorter and more  streamlined proof of the main step of \cite{brahma-pared}. This is the content of Theorem \ref{main}.\\
3) To generalize a reduction Theorem of Klarreich \cite{klarreich} to the context of relative hyperbolicity. This is the content of Theorem \ref{red}.\\

\smallskip

The main tool, Theorem \ref{red}, is a `reduction Theorem' ((3) above) which  allows us to deduce the existence and structure of Cannon-Thurston maps for the inclusion of one
 relatively hyperbolic metric space  into another, once we know the existence and structure of Cannon-Thurston maps 
for inclusions of certain relatively quasiconvex subspaces into {\it ends}. The exact statement of  
Theorem \ref{red} is somewhat technical. Suffice to say, this
is the appropriate relative hyperbolic generalization of inclusions of geometrically finite hyperbolic 3-manifolds  $M_{gf}$ into degenerate 
hyperbolic 3-manifolds $N^h$ such that \\
a) the inclusion of a boundary component $S_{gf}$ of $M_{gf}$ into the end $E^h$ of $N^h$ it bounds is a homotopy equivalence. \\
b) Each $S_{gf}$ is incompressible in $M_{gf}$.\\

\smallskip

We give the main
application below.

\noindent { \bf Theorem \ref{main2}: }{\it 
Suppose that $N^h \in H(M,P)$ is a hyperbolic structure 
on a pared manifold $(M,P)$ with incompressible boundary $\partial_0 M$. Suppose further that $N^h$ is not doubly degenerate. Let
$M_{gf}$ denotes a geometrically finite hyperbolic structure adapted
to $(M,P)$. Then the map  $i: \widetilde{M_{gf}}
\rightarrow \widetilde{N^h}$ extends continuously to the boundary
$\hat{i}: \widehat{M_{gf}}
\rightarrow \widehat{N^h}$.

 Let  $\partial {i}:  \partial {\til{M_{gf}}}
\rightarrow \partial {\til{N^h}}$ be the resulting Cannon-Thurston map extending  $i: \widetilde{M_{gf}}
\rightarrow \widetilde{N^h}$. 
Then $\partial i (a) = \partial i (b)$ for $a \neq b$ if and only if $(a,b) \in \RR$, where $\RR$ is the smallest 
closed equivalence relation containing the equivalence relations  generated by lifts of the ending laminations to 
$\til{M_{gf}}$.}

\smallskip

The last statement is informally abbreviated by saying that the Cannon-Thurston map identifies precisely the end-points of
leaves of the ending laminations. (Note that we have to pass to the transitive closure to get a precise statement.) It is curious that the doubly degenerate
case (dealt with in \cite{mahan-split, mahan-elct}) is the single exceptional case not amenable to the techniques of this paper.

\medskip

The last step of the programme of proving the existence of Cannon-Thurston maps for arbitrary finitely generated Kleinian groups and describing their structure
is dealt with in \cite{mahan-kl}.

\smallskip

\noindent { \bf Acknowledgments: } The authors would like to thank the referee for a careful reading and detailed and helpful comments and also for pointing out a
gap in an earlier draft.

\section{Background}

\subsection{Relative Hyperbolicity and Quasiconvexity}

Let $(X,d)$ be a path metric space. A collection of closed
 subsets $\HH = \{ H_\alpha\}$ of $X$ will be said to be {\bf uniformly
 separated} if there exists $D > 0$ such that
$d(H_1, H_2) \geq D$ for all distinct $H_1, H_2 \in \HH$.

\begin{defn} (Farb \cite{farb-relhyp})
The {\bf electric space} (or coned-off space) $\EE {(X, \HH )}$
corresponding to the
pair $(X,\HH )$ is a metric space which consists of $X$ and a
collection of vertices  $v_\alpha$ (one for each $H_\alpha \in \HH$)
such that each point of $H_\alpha$ is joined to (coned off at)
$v_\alpha$ by an edge of length $\half$. The sets $H_\alpha$ shall be
referred to as  horosphere-like sets and the vertices $v_{\alpha}$ as cone-points. \\
$X$ is said to be {\bf weakly hyperbolic} relative to the collection $\HH$ if $\EE {(X, \HH )}$
is a hyperbolic metric space.
\label{el-space}
\end{defn}

\begin{defn}
 A path $\gamma$  is said to be an electric geodesic (resp. electric $K$-quasigeodesic) if it is a geodesic (resp. $K$-quasigeodesic) in $\EXH$.\\
 $\gamma$ is said to be an  electric
$K$-quasigeodesic in (the electric space) $\EE {(X, \HH )}$
 {\bf without  backtracking}  if
 $\gamma$ is an electric $K$-quasigeodesic in $\EE {(X, \HH )}$ and
 $\gamma$ does not return to  any {\em horosphere-like set} $H_\alpha$ after leaving it.
\end{defn}

Let $i: X \rightarrow \EE {(X, \HH )}$ denotes the natural inclusion of spaces. Then
for a path $\gamma \subset X$, the path $i ( \gamma )$ lies in $\EE {(X, \HH )}$. Replacing maximal subsegments
$[a,b]$ of $i (\gamma )$  lying in a particular $H_\alpha$ by a path that goes from
$a$ to $v_\alpha$ and then from $v_\alpha$ to $b$, and repeating this for every $H_\alpha$ that
$i ( \gamma )$ meets we obtain a new path $\hat{\gamma}$. If $\hat{\gamma}$ is an electric geodesic (resp. $P$-quasigeodesic), $\gamma$ is called a {\em relative geodesic} (resp.
{\em relative $P$-quasigeodesic}). Paths (resp. geodesics or quasigeodesics) in $X$ shall be referred to as ambient paths
 (resp. geodesics or quasigeodesics). As above, an ambient path is said to be   without  backtracking if it
does not return to  any  horosphere-like set $H_\alpha$ after leaving it.
We shall usually be concerned with the case that $\gamma$ is an ambient geodesic/quasigeodesic without backtracking.

\begin{defn}
Relative  $P$-quasigeodesics in
$(X,\HH )$ are said to satisfy {\bf bounded region penetration} if there exists $B = B(P )$ so that for any two
    relative  $P$-quasigeodesics without backtracking
$\beta$, $\gamma$,
   joining $x, y$ we have
  \\
{\bf Similar Intersection Patterns 1:}  if
  precisely one of $ \beta , \gamma $ meets
 a  horosphere-like set $H_\alpha$,
then the length of this path (measured in the intrinsic path-metric
  on  $H_\alpha$ ) from the first (entry) point
  to the last
  (exit) point (of the relevant path) is at most $B$. \\
 {\bf Similar Intersection Patterns 2:}  if
 both $\beta , \gamma$ meet some  $H_\alpha $
 then the distance (measured in the intrinsic path-metric
  on  $H_\alpha$ ) from the entry point of
 $\beta$ to that of $\gamma$ is at most $B$; similarly for exit points. \\
\end{defn}

Replacing `$P$-quasigeodesic' by `geodesic' in the above definition, we obtain the notion of
relative  geodesics in
$(X,\HH )$   satisfying  bounded region penetration.

Families of paths which enjoy the above properties shall be said to have similar intersection patterns with horospheres.

\begin{defn} (Farb \cite{farb-relhyp} ) $X$ is said to be hyperbolic relative to the uniformly separated collection $\HH$ if \\
1) $X$ is weakly hyperbolic relative to $\HH$. \\
2) For all $P \geq 1$, relative $P$-quasigeodesics without backtracking satisfy the bounded penetration property. \\
Elements of $\HH$ will be referred to as {\bf horosphere-like} sets.
\end{defn}

\noindent {\bf Gromov's definition of  relative hyperbolicity \cite{gromov-hypgps} :}
 \begin{defn} ({\bf Gromov})
For any geodesic metric space
$(H,d)$, the {\em hyperbolic cone} (analog of a horoball)
$H^h$ is the metric space
$H\times [0,\infty) = H^h$ equipped with the
path metric $d_h$ defined by \\
1) $d_{h,t}((x,t),(y,t)) = 2^{-t}d_H(x,y)$, where $d_{h,t}$ is the induced path
metric on $H\times \{t\}$.  Paths joining
$(x,t),(y,t)$ and lying on  $H\times \{t\}$
are called {\em horizontal paths}. \\
2) $d_h((x,t),(x,s))=\vert t-s \vert$ for all $x\in H$ and for all $t,s\in [0,\infty)$, and the corresponding paths are called
{\em vertical paths}. \\
3)  For all $x,y \in H^h$,  $d_h(x,y)$ is the path metric induced by the collection of horizontal and vertical paths. \\
\end{defn}

\begin{defn}
Let $X$ be a geodesic metric space and $\HH$ be a collection of mutually disjoint uniformly separated subsets of $X$. The space
$X$ is said to be  hyperbolic relative to $\HH$ in the sense of Gromov, if the quotient space $\GG (X, \HH)$,  obtained by attaching the hyperbolic cones
$ H^h$ to $H \in \HH$  by identifying $(z,0)$ with $z$
for all $H\in \HH$ and $z \in H$,
 is a complete hyperbolic metric space. The collection $\{ H^h : H \in \HH \}$ is denoted
as ${\HH}^h$. The induced path metric is denoted as $d_h$.

We shall refer to $\GG (X, \HH)$ as the {\bf Gromov cone} for the pair $(X, \HH)$.
\end{defn}

\begin{theorem} (Bowditch \cite{bowditch-relhyp})\label{bow-rel}
The following are equivalent: \\
1) $X$ is  hyperbolic relative to the
collection $\HH$ of uniformly separated subsets of $X$.  \\
2) $X$ is  hyperbolic relative to the
collection $\HH$ of uniformly separated subsets of $X$ in the sense of Gromov. \\
3) $\GXH$ is  hyperbolic relative to the
collection $\HH^h$. \\
\end{theorem}

\begin{defn} Let $X$ be hyperbolic relative to the collection $\HH$.
We call a set $W \subset X$  relatively $K$-quasiconvex 
 if \\
1) $W$ is hyperbolic relative to the collection $\WW = \{ W \cap H: H \in \HH \}$ \\
2) $\EE (W, \WW )$ is $K$-quasiconvex in $\EXH$.  \\
$W \subset X$ is  relatively quasiconvex if it is relatively $K$-quasiconvex for some $K$.
\end{defn}

\noindent {\bf Ends:}\\
Let $Y$ be hyperbolic rel. $\HH$.
Now let $\mathcal{B}=\{B_\alpha\},\alpha\in\varLambda$, for some indexing set $\varLambda$, 
be a collection of uniformly relatively quasiconvex sets inside $Y$. 
Here each $B_\alpha$ is relatively quasiconvex with respect to 
the collection $\{B_{\alpha\beta}\}$, given by $B_{\alpha \beta} = B_\alpha \cap H_\beta$. 
We also assume that the sets $H_\beta$ are  uniformly $D$-separated.\\

\begin{defn} Let $Y$ be hyperbolic relative to the collection $\HH$ and $X$ be strongly hyperbolic with respect to a collection $\JJ$.
A map $i: Y \rightarrow X$ is said to be strictly type-preserving if: \\
1) For every $H \in \HH$, $i(H) \subset J_H$ for some $J_H\in \JJ$. \\
2) For every $J \in \JJ$, $i^{-1} (J) = \emptyset$ or $i^{-1} (J) = H_J$ for some $H_J \in \HH$
\end{defn}

A map of path-metric spaces is a length-isometry if it preserves lengths of paths.
\begin{defn}\label{ei}
A strictly type-preserving length-isometric inclusion $i: Y \hookrightarrow X$   of relatively hyperbolic metric spaces is said to be an {\bf ends-inclusion} if \\
1) There exist collections  $\JJ = \{ J_\beta \}$, $\HH = \{ H_\beta \}$ such that $X$ is hyperbolic rel. $\JJ$ and $Y$ is hyperbolic rel. $\HH$
(note that $\beta$ ranges in the same indexing set).\\
2) There exists a collection $\BB =\{ B_\alpha \},\alpha\in\varLambda$,   of relatively quasiconvex subsets of $Y$. 
Each $B_\alpha$ is relatively quasiconvex with respect to 
the collection $\{B_{\alpha\beta}\}$ given by $B_{\alpha \beta} = B_\alpha \cap H_\beta$. \\
3) There exists a collection    $\mathcal{F}=\{F_\alpha \subset X\},\alpha\in\varLambda$, of relatively quasiconvex subsets of $X$
(thought of as {\bf ends} of $X$),
such that
$B_\alpha = F_\alpha \cap Y$, $\forall \alpha$ and  $X = Y \cup 
\{\bigcup_\alpha F_\alpha\}$. We also have the inclusion maps  $i_\alpha : B_\alpha \rightarrow F_\alpha$. \\
4) Each $F_\alpha$ is strongly hyperbolic relative to the collection $\{ F_{\alpha \beta} = F_\alpha \cap J_\beta\}$. \\
5) If $\HH_0$ is the subcollection of elements $H_\gamma \in \HH$ such that $H_\gamma \cap F_\alpha =\emptyset$ for all $F_\alpha$, then
$\JJ = \HH_0 \cup\bigcup_{\alpha , \beta} \{ F_{\alpha \beta} \}$. \\ We let $\HH_1 = \HH \setminus \HH_0$.
 \end{defn}

\begin{rmk}\label{mot-eg}{\rm 
It might be useful here to keep the motivating example of a pared hyperbolic 3-manifold $N$ with incompressible boundary (cf. Definition
\ref{pared} below) in mind. We give an informal sketch of the setup to fix notions.
In this situation, there exists a geometrically finite manifold $M$ and an embedding $i: M \rightarrow N$
such that $N \setminus M$ consists of finitely many products of the form $S \times [0, \infty)$ for $S$ a finite area hyperbolic surface.
Then $X$ (resp. $Y$) in Definition \ref{ei} corresponds to the universal cover of $N$ (resp. $M$).
The lifts of the  $S \times \{ 0 \}$'s  correspond to $\{ B_{\alpha} \}$.
The lifts of the cusps of the $S \times [0, \infty)$'s  correspond to $\{ F_{\alpha \beta} \}$.
There might be cusps in $M$ which have no curves parallel to  the cusps of the  $S \times \{ 0 \}$'s.
Lifts of such cusps correspond to  $\HH_0$. Finally, the lifts of the cusps of the  $S \times \{ 0 \}$'s  correspond to $\HH_1$.

$M$ is often referred to as the {\bf relative Scott core} of $N$.

}
\end{rmk}

\begin{rmk}\label{geods-ei}{\rm Note that the ends-inclusion $i: Y \hookrightarrow X$ induces an {\bf isometric embedding}
$\hat{i}: \EE (Y, \BB) \rightarrow \EE (X, \FF)$. Further, every point of $\EE (X, \FF)$ is within bounded 
distance (in fact distance $\half$) of the image of  $\EE (Y, \BB)$. The points of $\EE (X, \FF)$} not {\rm in the 
image of  $\EE (Y, \BB)$ correspond precisely to points of $F_\alpha \setminus i_\alpha ( B_\alpha)$ for some $\alpha$.
It follows that for any electric geodesic (resp. $P-$ quasigeodesic) $\sigma$ in $\EE (Y, \BB)$, $\hat{i} (\sigma)$ is an 
 electric geodesic (resp. $P-$ quasigeodesic) in $\EE (X, \FF)$. } \end{rmk}

\begin{lemma} \cite{bowditch-relhyp}\label{eg2}  Let  $X$ be  a hyperbolic metric space and  let $\BB$ be a collection of  uniformly separated
uniformly quasiconvex sets.  Then $X $ is weakly hyperbolic relative to the collection $\mathcal{B}$. 
\end{lemma}

Let $X$ be a $\delta$-hyperbolic metric
space, and $\mathcal{B}$ a family of $C$-quasiconvex, $D$-separated,
 collection of subsets. Then by Lemma \ref{eg2} (see also \cite{farb-relhyp}),
$X_{el} = \EE(X,\BB )$ obtained by electrocuting the subsets in $\mathcal{B}$ is
a $\Delta = \Delta ( \delta , C, D)$ -hyperbolic metric space. Now,
let $\alpha = [a,b]$ be a hyperbolic geodesic in $X$ and $\beta $ be
an electric 
$P$-quasigeodesic without backtracking
 joining $a, b$. Replace each maximal subsegment, (with end-points $p,
 q$, say)
starting from the left
 of
$\beta$ lying within some $H \in \mathcal{H}$
by a hyperbolic  geodesic $[p,q]$. The resulting
{\bf connected}
path $\beta_{ea}$ is called an {\em electro-ambient path representative} in
$X$. 

\smallskip

Note that $\beta_{ea}$ {\it need
  not be a hyperbolic quasigeodesic}. However, the proof of Proposition
  4.3 of Klarreich \cite{klarreich} gives the following. (See  \cite[Lemma 2.5]{mahan-split} for a proof
of the forward direction. 
The converse direction follows directly from the proof of   \cite[Proposition
  4.3]{klarreich}.) 

\begin{lemma}  
Given $\delta$, $C, D, P$ there exists $C_3$ such that the following
holds: \\
Let $(X,d)$ be a $\delta$-hyperbolic metric space and $\mathcal{H}$ a
family of $C$-quasiconvex, $D$-separated collection of quasiconvex
subsets. Let $(X,d_e)$ denote the electric space obtained by
electrocuting elements of $\mathcal{H}$.  Then, if $\alpha , \beta_{ea}$
denote respectively a hyperbolic geodesic and an electro-ambient
$P$-quasigeodesic with the same end-points, then $\alpha$ lies in a
(hyperbolic) 
$C_3$ neighborhood of $\beta_{ea}$.

Conversely, given a hyperbolic geodesic $\alpha$, there exists  an electro-ambient $P-$ quasigeodesic $\gamma_{ea}$
lying in  a
(hyperbolic) 
$C_3$ neighborhood of $\alpha$.

\label{ea-strong}
\end{lemma}

We shall abbreviate this  as:\\
{\it  Hyperbolic geodesics lies hyperbolically close to electro-ambient representatives of electric geodesics joining their end-points.
Conversely, given a hyperbolic geodesic there is an electro-ambient quasigeodesic lying close to it.}

\smallskip

A word of clarification here regarding the hypotheses of Lemma \ref{ea-strong}. $D$- separatedness is only a technical assumption.
Given $X, \mathcal{H}$, let  $X_1 = X \bigcup_{H \in \mathcal{H}} (H \times [0,1])$, equipped with the quotient topology, where $(h,0) \in  (H \times [0,1])$
is identified with $h \in H \subset X$. Then the collection $\{ H \times \{ 1\} : H \in \mathcal{H}\}$ is automatically 2-separated and the inclusion of $X$
in $Y$ is a quasi-isometry. However, the requirement that each $H$ is $C-$ quasiconvex is an essential assumption and the conclusion of Lemma \ref{ea-strong}
fails without this assumption. It is {\em not} sufficient to assume that $X$ is (weakly) hyperbolic relative to the collection  $\mathcal{H}$.
A simple counterexample is given by a doubly degenerate 3-manifold, with the 2 ends corresponding to the 2 elements of $\mathcal{H}$. We are grateful to the referee
for bringing this to our notice.

\subsection{Cannon-Thurston Maps}

For a hyperbolic metric space $X$, the Gromov bordification will be denoted by $\bbar X$.

\begin{defn}
Let $X$ and $Y$ be hyperbolic metric spaces and
$i : Y \rightarrow X$ be an embedding.
 A {\bf Cannon-Thurston map} $\bbar{i}$  from $\bbar{Y}$ to
 $\bbar{X}$ is a continuous extension of $i$ to the Gromov bordifications $\bbar{X}$ and $\bbar{Y}$. \end{defn}

 The following lemma from  \cite{mitra-trees} gives a necessary and sufficient condition for the existence of Cannon-Thurston maps.

\begin{lemma}\cite{mitra-trees}\label{ct-crit}
 A Cannon-Thurston map $\bbar{i}$ from $\bbar Y $ to $\bbar{X}$ exists for the proper embedding $i\colon Y\to X$ if and only if
there exists a non-negative function $M(N)$ with $M(N)\rightarrow \infty$ as $N\rightarrow \infty$ such that the following holds: \\
Given $y_0 \in Y$, for all geodesic segments $\lambda$ in
$Y$ lying outside an $N$-ball around $y_0$ $\in Y$, any geodesic segment in $X$ joining the end points of $i(\lambda )$ lies outside the
$M(N)$-ball around $i(y_0)\in X$.
\end{lemma}

Note that due to stability of quasigeodesics, the above statement is also true if geodesics are replaced by uniform quasigeodesics.\\

Let $X$ and $Y$ be    hyperbolic relative to the collections $\HH_X$ and $\HH_Y$ respectively. Let 
$\widehat{X} = \EE (X, \HH_X), \widehat{Y} = \EE (Y, \HH_Y)$.
Let $i\colon Y\to X$ be a strictly type-preserving proper embedding. 
Then
the proper embedding $i\colon Y\to X$  induces a proper embedding  $i_h\colon \GG(Y,\HH_Y)\to \GG(X,\HH_X)$ and a map $\hat{i}: \widehat{X}
\rightarrow \widehat{Y}$.

\begin{defn}
A Cannon-Thurston map is said to exist for the
 pair $X, Y$ of relatively hyperbolic metric spaces and a strictly type-preserving inclusion
$i: Y \rightarrow X$ if a Cannon-Thurston map  exists for the induced map $i_h\colon \GG(Y,\HH_Y)\to \GG(X,\HH_X)$.
\end{defn}

In \cite{mj-pal} Lemma \ref{ct-crit} was generalized to relatively hyperbolic metric spaces as follows.

\begin{lemma}\label{crit-relhyp} (\cite{mj-pal} Lemma 1.28)
Let $Y, X$ be hyperbolic rel. $\YY , \XX$ respectively. Let $Y^h = \GG (Y, \YY ),
\widehat Y = \EE (Y, \YY )$ and $X^h = \GG (X, \XX ),
\widehat X = \EE (X, \XX )$.
A Cannon-Thurston map for $i\colon Y \to X$ exists if and only if
there exists a non-negative function $M(N)$ with
$M(N)\rightarrow \infty$ as $N\rightarrow \infty$ such that the
following holds: \\
Suppose $y_0\in Y$,  and $\hat \lambda$ in $\widehat Y$ is
an electric geodesic segment  starting and ending  outside horospheres.
If  $\lambda^b = \hat \lambda \setminus \bigcup_{K \in \YY} K$
 lies outside  $B_N (y_0) \subset Y$,
then for any electric quasigeodesic $\hat \beta$ joining the
end points of $\hat i (\hat \lambda)$ in $\widehat X$,
$\beta ^b = \hat \beta \setminus \bigcup_{H \in \XX} H$ lies outside
 $B_{M(N)} (i(y_0)) \subset X$. 
\end{lemma}    

The above necessary and sufficient condition for existence of Cannon-Thurston map for relatively 
hyperbolic spaces can also be used as a definition of Cannon-Thurston map for relatively hyperbolic spaces.
Hence the following definition makes sense.

\begin{defn}
 A collection of proper, strictly type preserving embedding $i_\alpha:Y_\alpha \rightarrow X_\alpha$ 
of relatively hyperbolic spaces is said to extend to a collection of
{\bf uniform Cannon-Thurston maps} if there exists $M(N) \rightarrow \infty$ as $N\rightarrow \infty$
such that the functions $M_\alpha(N)$ (obtained in Lemma \ref{crit-relhyp} above) satisfy $M_\alpha(N) \geq M(N)$ for all $\alpha$.
\end{defn}

We shall often abbreviate {\bf Cannon-Thurston} as {\bf CT} in what follows.
Lemma \ref{crit-relhyp} says that it is enough to consider only the \textquoteleft bounded\textquoteright 
-part of the electric quasigeodesic in a relatively hyperbolic space $X$ in order to prove existence 
of \textbf{CT} map. For ease of reference below, we make the following definition.

\begin{defn}
Let $ X$ be hyperbolic rel. $\XX$. If $\sigma$ is a path in $X$, the {\bf bounded part} $\sigma^b$ of $\sigma$
with respect to $(X, \XX)$
is defined as $\sigma \setminus \bigcup_{H \in \XX} H$.
\end{defn}

If there is no ambiguity, we shall  refer to the   bounded part of $\sigma$
with respect to $(X, \XX)$  simply as the  bounded part of $\sigma$.

We shall use the notion of electro-ambient path representatives to obtain an alternate criterion for the
existence of Cannon-Thurston maps in the case of an {\bf ends-inclusion}. Combining Lemma \ref{crit-relhyp} with 
Lemma \ref{ea-strong} we have the following.

\begin{lemma} \label{finalcrit} Let $X, Y$ be hyperbolic rel. $\JJ , \HH$ respectively and
$i: Y \rightarrow X$ be an ends-inclusion of relatively hyperbolic spaces.
A Cannon-Thurston map for
$i\colon Y \hookrightarrow X$ 
exists if and only  there  exists  a  non-negative  function  
$M(N)$ with $M(N)\rightarrow \infty$ as $N\rightarrow \infty$ such that the following holds. \\
Suppose $y\in Y$,  and $\hat \lambda$ in $\widehat Y$ is
an electric geodesic segment  starting and ending  outside horospheres,such that $\lambda^b = \hat \lambda 
\setminus \bigcup_{K \in \HH} K$, the bounded part of $\hat \lambda$ lies outside  $B_N (y) \subset Y$.

Then for some electric quasigeodesic $\hat \rho$ joining the end points 
of $\hat i (\hat \lambda)$ in $\widehat X$, the bounded part $\rho^b_{ea} = \hat{\rho}_{ea} \setminus \cup _
{H\in\JJ}H$ of the electro-ambient representative $\rho_{ea}$ (of $\hat \rho$) lies outside $B_{M(N)} (i(y)) \subset X$.
\end{lemma}

\subsection{Pared Manifolds}\label{pared-sec}
The main examples of interest in this paper are pared 3-manifolds. 

\begin{defn} \label{pared} A {\bf pared manifold} is a pair $(M,P)$, where $M$ is a compact  irreducible 3-manifold with boundary
$\delta M$ and $P
\subset \delta M$ 
is a (possibly empty) 2-dimensional submanifold  with boundary (of $\delta M$) such that \\
\begin{enumerate}
\item Any $\pi_1$-injective map of a torus or Klein bottle into $M$ is homotopic to a map into $\delta M$.
\item The fundamental group of each component of $P$ injects into the
fundamental group of $M$.
\item The fundamental group of each component of $P$ contains an abelian 
subgroup with finite index.
\item Any map $C: (S^1 \times I, \delta (S^1 \times I)) \rightarrow (M,P)$
such that $\pi_1 (C)$ is injective, is homotopic {\it rel} boundary to $P$.
\item $P$ contains every component of $\delta M$ which has an abelian subgroup
of finite index.
\end{enumerate}
\end{defn}

 A pared manifold $(M,P)$ is said to have {\bf
  incompressible boundary} 
if each component of $\delta_0 M = \delta M \setminus P$ is
incompressible in $M$.

Further, $(M,P)$ is said to have {\bf
 no accidental parabolics} if 
\begin{enumerate}
\item It has incompressible boundary.
\item  If some curve $\sigma$ on a component of $\delta_0 M$ is freely 
homotopic in $M$ to a curve $\alpha$ on a component of $P$, then
  $\sigma$ is homotopic to $\alpha$ in $\delta M$.  
\end{enumerate}

\begin{defn} \cite{thurston-hypstr1, thurston-hypstr3} A hyperbolic structure on a pared manifold $(M,P)$ is defined to be a complete hyperbolic
structure 
on the interior of $M$ given by a discrete faithful representation
$\rho: \pi_1(M) \rightarrow Isom ({\mathbb{H}}^3)$ such that any homotopically nontrivial loop in $M$ represented by a parabolic
is homotopic into $P$. Further, for any component $P_i$ of $P$, and  any homotopically essential curve $\gamma$ in $\pi_1(P_i)$ ($\subset \pi_1(M)$),
$\rho (\gamma)$ is a parabolic.
\\
The space of hyperbolic structures on $(M,P)$ is denoted by $H(M,P)$. \end{defn}

Let $\Gamma = \rho ( \pi_1(M)) \subset Isom ({\mathbb{H}}^3)$. A hyperbolic structure on $(M,P)$ is said to be geometrically
finite (resp. infinite) if $\Gamma$ is a geometrically
finite (resp. infinite) Kleinian group. Thurston's hyperbolization theorem \cite{thurston-hypstr1, thurston-hypstr3, kapovich-book,otal-geom} ensures that
$H(M,P)$ contains a geometrically finite structure $M_{gf}$.  Further, the limit set of a geometrically finite $\Gamma$ is equivariantly homeomorphic
to the boundary of the Gromov cone $\GXH$ where $X$ is the universal cover $\til{M}$ and the parabolic subgroups $\HH$ correspond 
to the fundamental groups
of the components of $P$. 
Very often, in what follows we shall not be considering all of ${\mathbb{H}}^3/\Gamma$ but rather its convex core, or equivalently, the
quotient of the convex hull of the limit set of $\Gamma$ by $\Gamma$. By slight abuse of notation, we shall continue
to denote the convex core of 
$M_{gf}$ by  
$M_{gf}$. 

We give a slightly different but equivalent description of  accidental parabolics in terms of hyperbolic structures on $(M,P)$. Recall (Definition
\ref{pared} above) that for a pared manifold $(M,P)$, any map $C: (S^1 \times I, \delta (S^1 \times I)) \rightarrow (M,P)$
of an annulus
such that $\pi_1 (C)$ is injective, is homotopic {\it rel} boundary to $P$.
An element $\gamma \in \Gamma$ is an accidental parabolic, if the converse is false, i.e.\\
a) If there exists a homotopically essential
map $C: (S^1 \times I, \delta (S^1 \times I)) \rightarrow (M,P)$  such that $C(S^1 \times \{0 \})$ is contained in $\delta_0 M$, 
$C(S^1 \times \{1 \})$ is contained in a component $P_i$ of $P$, but $C$ is {\em not} homotopic rel. boundary to a map 
$ (S^1 \times I, \delta (S^1 \times I)) \rightarrow \delta M$. \\
b) a geodesic representative of $\gamma$ in $M$ is freely homotopic to the core curve of the annulus 
$C (S^1 \times I)$.

A component $P_i$ of $P$ for which such a map $C$ exists is called {\it exceptional}.

In summary an accidental parabolic is given by the core curve of a homotopically essential
map $C: (S^1 \times I, \delta (S^1 \times I)) \rightarrow (M,P)$ of an annulus into a pared manifold
$(M,P)$  such that \\
a) $C(S^1 \times \{0 \})$ is contained in $\delta_0 M (= \delta M \setminus P)$, \\
b) $C(S^1 \times \{1 \})$ is contained in a component $P_i$ of $P$, \\
c) $C$ is {\em not} homotopic rel. boundary to a map 
$ (S^1 \times I, \delta (S^1 \times I)) \rightarrow \delta M$.

For a hyperbolic structure $N^h \in H(M,P)$ adapted to $(M,P)$, an {\bf exceptional cusp} is a cusp corresponding to an exceptional component $P_i$.
{\bf Exceptional horoballs} are lifts of (neighborhoods of) exceptional cusps. Boundaries of exceptional horoballs
are called {\bf exceptional horospheres}.

We now describe how to adjoin exceptional cusps to ends having accidental parabolics so that the resulting set can be treated on an equal footing with
ends containing no accidental parabolics.
 
Let $E$ be an end of $N^h$ and $\Sigma \subset \delta M$ be its boundary. Let $\sigma_1, \cdots \sigma_k \subset \Sigma$ be all the simple closed
curves on $\Sigma$ corresponding to accidental parabolics. Then each $\sigma_i$ is homotopic into an  exceptional cusp and there is an embedded annulus $A_i$
with one boundary component $\sigma_i$ and the other component $\sigma_i^\prime$ in the exceptional cusp.
We choose $\sigma_i^\prime$  to be geodesic in the canonical 
flat metric on the
boundary of the exceptional cusp. Then  $\sigma_i^\prime$  bounds a totally geodesic annulus $C_i$ contained  in the exceptional
cusp bounded by  $\sigma_i^\prime$ and  isometric to the quotient of a
2-dimensional horodisk by a cyclic parabolic group. Note that if the  exceptional cusp is rank one, then $C_i$ equals the   exceptional cusp. The union
$E \bigcup_i (A_i \bigcup_i C_i)$ will be termed an {\bf augmented end}.

We shall need 
the Thurston-Canary covering theorem \cite{Thurstonnotes}[Ch. 9]
\cite{canary-cover} in the context of pared manifolds. The version below combines the covering theorem with the tameness theorems of
Bonahon, Agol and Gabai-Calegari  \cite{bonahon-bouts,agol-tameness,gab-cal}.

\begin{theorem}  \cite{Thurstonnotes}
\cite{canary-cover} Let $M = {\Hyp}^3/\Gamma$ 
be a complete hyperbolic 3-manifold. A finitely generated subgroup
${\Gamma}'$ is geometrically infinite if and only if it  contains
a finite index subgroup of a geometrically infinite peripheral
subgroup.
\label{cover}
\end{theorem}

Another fact we shall need in this context is the following (see also \cite{canary-cover}):

\begin{lemma} 
 \label{canary-qc}  Let $E$ be an augmented degenerate end for a hyperbolic structure $N^h$ on a pared manifold $(M,P)$ with
incompressible boundary. Let $\til E$ be a lift of $E$
to $\til M$, equipped with this  hyperbolic structure. Then $\til E$ is not relatively quasiconvex in $\til{N^h}$ if and only if there is a 
component $F$  of $\delta_0 M$ such that \\
\begin{enumerate}
 \item $F$ bounds a degenerate end other than $E$, and
\item $F$ is homotopic into $E$.
\end{enumerate}
\end{lemma}

\begin{proof}  The proof is essentially a rerun of some of the arguments appearing in \cite{bonahon-bouts}. 
 Suppose $\til E$ is not relatively quasiconvex in $\til{N^h}$. Then there exists a sequence of  closed curves $\sigma_i$ on $\partial E$ whose geodesic realizations $\gamma_i$
in $N^h$ satisfy $d(\gamma_i, E) \rightarrow \infty$ as $i \rightarrow \infty$ (Section 2.2 of  \cite{bonahon-bouts} proves the existence of closed curves satisfying the above
property). Then (cf. Proposition 5.1 of \cite{bonahon-bouts}) a subsequence of the $\sigma_i$'s converges to an ending lamination $\Lambda$ on $\partial E$.
If $\Lambda_E$ is the ending lamination for the end $E$, then $\Lambda$ is different from $\Lambda_E$.  

If the support of $\Lambda$ is all of $\partial E$ and $\Lambda$ contains no simple closed curve, then $N^h$ is doubly degenerate.
Else any simple closed curve in $\Lambda$ gives rise to an accidental parabolic. Let $F$ be a  connected subsurface of $\partial E$ supporting an ending lamination contained in $\Lambda$.
Then $F$ satisfies the conclusions of the Lemma.  \end{proof}

\section{Reduction Theorem and Kleinian Groups}

\subsection{The Main Theorem}
Before stating the main Reduction Theorem \ref{red} below, we briefly sketch the proof idea in the special case of
 hyperbolic 3-manifolds $N$
 with incompressible core $M$ and no parabolics. For concreteness, suppose that
$N$  has one end and that the end $E= N \setminus M$ is homeomorphic to $S \times [0,\infty)$
for a compact hyperbolic surface $S$. Theorem \ref{red} says in this case that if the inclusion of $\til{S}$ into 
$\til E$ has a CT map, then so does the inclusion of $\til{M}$ into 
$\til N$. 
The proof idea is as follows:\\ Let $\EE$ denote the collection of lifts of the end $E$ to $\til N$ and let $\SSS$ denote the collection of lifts of $S$ to $\til M$.
Then by Lemma \ref{canary-qc}, each lift $E_\alpha \in \EE$ is relatively quasiconvex in $\til N$.

Let $[a,b] \subset \til{M}$ be a geodesic in the intrinsic path metric on $\til M$ lying outside a large ball about a fixed reference point 
$m \in \til{M}$.
We want to construct an electro-ambient $P-$quasigeodesic with respect to $(\til{N}, \EE)$ lying outside a large ball in $\til N$. Towards this, first construct 
an  electro-ambient $P-$quasigeodesic $[a,b]_q$ with respect to $(\til{M}, \SSS)$ in $\til M$ joining $a, b$ and
lying within a $P-$ neighborhood of $[a,b]$ (by
Lemma \ref{ea-strong}). This gives us a sequence of points $a=a_0, \cdots a_n = b$ such that the "odd subpaths"
$[a_{2i}, a_{2i+1}]_q$ of $[a,b]_q$ have interiors disjoint from the elements of $\SSS$, whereas the 
"even subpaths"
$[a_{2i+1}, a_{2i+2}]_q$ of $[a,b]_q$ lie entirely within some ${\til{S}}_\alpha \in \SSS$. Since all these are subpaths of 
$[a,b]_q$, they all lie outside a large ball about $m$. Now replace the even subpaths
$[a_{2i+1}, a_{2i+2}]_q$ by a geodesic $\bbar{a_{2i+1}, a_{2i+2}}$ in the corresponding ${\til{E}}_\alpha \in \EE$
joining $a_{2i+1}, a_{2i+2}$.  Since 
the inclusion of $\til{S}$ into 
$\til E$ has a CT map, it follows that each of the geodesic segments  $\bbar{a_{2i+1}, a_{2i+2}}$ lies outside a (uniformly) large ball
about $m$. Concatenate these together by interpolating the odd subpaths
$[a_{2i}, a_{2i+1}]_q$ of $[a,b]_q$. This gives an electro-ambient $P-$quasigeodesic $\bbar{a,b}^q$
with respect to $(\til{N}, \EE)$ by Remark \ref{geods-ei}. Further,  $\bbar{a,b}^q$ also lies outside a large ball
about $m$.
Finally, by Lemma \ref{ea-strong}, the hyperbolic geodesic in $\til N$ lies in a bounded neighborhood of  $\bbar{a,b}^q$ and hence 
lies outside a large ball
about $m$. Lemma \ref{ct-crit} now furnishes the required CT map. We now proceed with the general case.

\begin{theorem}\label{red}
 Let $Y,X, \HH, \JJ, \BB = \{B_\alpha\}, \FF = \{F_\alpha\}, {\FF}_\alpha = \{F_{\alpha \beta}\}$ be as in Definition \ref{ei}
and $i:Y \rightarrow X$ be an {\bf ends-inclusion} of spaces. 
Then the  ends inclusion $i:Y \rightarrow X$ extends to a Cannon-Thurston map if the inclusions $i_\alpha:B_\alpha \rightarrow F_\alpha$ 
extend {\bf uniformly} to  Cannon-Thurston maps for all $\alpha$.
\end{theorem}

\begin{proof} Fix a base point $y \in Y$ and consider a large enough 
ball $U_N(y) \subset Y$. Let $\hat{\eta} \subset \widehat{Y} = \EE(Y, \HH )$ be an electric geodesic segment,
starting and ending outside elements of $\HH$. Let $\eta^b$ denote 
 the bounded part of $\hat{\eta}$ with respect to $(Y, \HH)$, and assume that it lies 
outside $U_N(y) \subset Y$, \textit{i.e} $\eta^b\cap U_N(y) = \varnothing$.

Let $\BB_0 = \{B_\nu \in \BB : \overline{\eta^b} \cap B_\nu \neq \emptyset$ and $U_N(y) \cap B_\nu \neq \emptyset \} $, 
where $\overline{\eta^b}$ denotes the closure of ${\eta^b}$.
 
For each $B_\nu \in \BB_0$, let $\eta^b(\nu) = \eta^b \cap B_\nu$. Then $\eta^b(\nu)$ lies outside $U_N(y) \cap B_\nu$. 
Let $y_\nu$ be the nearest point projection of $y$ on $B_\nu$  in the metric $d_Y$
of $Y$. Since $F_\nu \cap Y = B_\nu$ it follows that 
$y_\nu$ is also (up to bounded error) a nearest point projection of $y$ on $F_\nu$  in the metric $d_X$
on $X$.
Then $y_\nu \in U_N(y) \cap B_\nu$. 
Let  $d_Y(y,y_\nu) = R_\nu$. Consider the ball $U_{(N-{R_\nu})}(y_\nu)$, of radius $N-R_\nu$ about $y_\nu$ in $Y$. 
 $U_{(N-{R_\nu})}(y_\nu) \cap B_\nu$ is a ball in $B_\nu (\subset Y)$ of radius $N-R_\nu$ based at 
$y_\nu$. We denote this ball as $U(\nu)$.  Then $\eta^b(\nu)\subset B_\nu \setminus U(\nu)$.

 Let $\hat \rho$ be the electric geodesic in $\widehat{X}= \EE(X, \JJ )$
 joining the end points of $\hat{i}(\hat \eta)$. Since $\widehat{Y}$ is weakly hyperbolic rel. $\BB$, it follows that
$\widehat{X}$ is weakly hyperbolic rel. $\FF$. Let the  electro-ambient 
path representative of $\hat \rho$ with respect to $\FF$ be $\hat \rho_{ea}$. Let $\rho^b_{ea}= \hat{\rho}_{ea} \setminus \cup _{J\in\JJ}J$
 be the bounded part of $\hat \rho_{ea}$ with respect to $(X, \JJ)$. By Remark \ref{geods-ei}, we may assume that $\rho^b_{ea} \setminus \cup
\bigcup_{F_\alpha \in \FF} F_\alpha = \eta^b  \setminus \cup
\bigcup_{B_\alpha \in \BB} B_\alpha$, i.e. $\rho^b_{ea}$ and $\eta^b$ coincide outside $\FF$. 

As per  hypothesis,
\textbf{CT} maps exist {\it uniformly} for each $B_\nu \hookrightarrow F_\nu$. By Lemma \ref{finalcrit}, there exists
a function $M(N) \rightarrow \infty$ as $N \rightarrow \infty$ such that $\rho^b_{ea}(\nu)$ lies outside $U_{M (N-R_\nu)}(y_\nu)$, 
$\forall \nu$. It is worth noting that the function $M(N)$ is independent of $\nu$ by definition of uniformity.

Since $Y$ is properly embedded in $X$, it follows that there exists a function 
$M_1(N) \rightarrow \infty$ as $N \rightarrow \infty$ such that if $ x, y \in Y$ and $d_Y(x,y) \geq N$ then $d_X(i(x),i(y)) \geq M_1(N)$.
It follows immediately that $\rho^b_{ea} \setminus \cup
\bigcup_{F_\alpha \in \FF} F_\alpha$ lies outside $U^X_{M_1(N)}(i(y))$.

Hence $\rho^b_{ea}(\nu)$ lies outside $U^X_{M_1(R_\nu) + M (N-R_\nu)}(i(y))$ -- a ball of radius $M_1(R_\nu) + M (N-R_\nu)$ in $X$,
i.e. $$d_X ( ( \rho^b_{ea}(\nu ) , i(y)) \geq M_1(R_\nu) + M (N-R_\nu)$$ for all $\nu$. 

Let $M_2(N) = {\rm inf}_\nu (M_1(R_\nu) + M (N-R_\nu))$, and $M_3(N) = min (M_1(N), M_2(N)$,
which is again a proper function of $N$, i.e. $M_3(N)\rightarrow \infty$ as $N \rightarrow \infty$.
 This proves that $\eta^b$ and $\rho^b_{ea}$ 
satisfy the criteria of Lemma \ref{finalcrit}.

Hence, the theorem follows. 
\end{proof}

An important fact  we used in the above proof is that $\rho^b_{ea} \setminus \cup
\bigcup_{F_\alpha \in \FF} F_\alpha = \eta^b  \setminus \cup
\bigcup_{B_\alpha \in \BB} B_\alpha$, i.e. $\rho^b_{ea}$ and $\eta^b$ 
can be chosen to coincide outside $\FF$.  This followed from Remark \ref{geods-ei}.

Now let $\partial i$ denote
the Cannon-Thurston map on the boundary $\partial \GG (Y, \HH )$ obtained in Theorem \ref{red}. We would like to know
exactly which points are identified by the CT map $\partial i$. Towards this, we set up some notation.

The inclusions $i_\alpha: B_\alpha \rightarrow F_\alpha$ induce CT maps $\partial i_\alpha : \partial \GG(B_\alpha, B_{\alpha \beta})
 \rightarrow  \partial \GG( F_\alpha, F_{\alpha \beta})$
 by the hypothesis of Theorem \ref{red}. Each such map $\partial i_\alpha$ induces an equivalence relation $\RR_\alpha$ on 
$\partial \GG(B_\alpha, B_{\alpha \beta})$ given
by $a \RR_\alpha b$ if and only if $\partial i_\alpha (a) = \partial i_\alpha (b)$. Since $\GG(B_\alpha, B_{\alpha \beta})$
 is quasiconvex in $\GG (Y, \HH )$ it follows that
$\partial \GG(B_\alpha, B_{\alpha \beta})$ embeds homeomorphically in $\partial \GG (Y, \HH )$. Hence  $\RR_\alpha$ induces a natural equivalence relation
(also denoted as $\RR_\alpha$) on $\partial \GG (Y, \HH )$ by identifying points on $\partial \GG(B_\alpha, B_{\alpha \beta})$
 with their images under inclusion in 
$\partial  \GG (Y, \HH )$. We shall call the relation $\RR_\alpha$  on $\partial  \GG (Y, \HH )$ the {\bf CT relation induced by $i_\alpha$}. Let $\RR_t$
denote the transitive closure of the union $\bigcup_\alpha \RR_\alpha$. Finally, let $\RR$ denote the closure of $\RR_t$ thought of as a subset
of $\partial  \GG (Y, \HH ) \times \partial  \GG (Y, \HH )$ with the product topology. Thus $\RR$ is the smallest closed equivalence relation generated by
the $\RR_\alpha$'s. 

As in the discussion preceding Theorem \ref{red}, we give a quick sketch of what goes on in the special case of a
hyperbolic 3-manifold $N$
 with incompressible core $M$, one simply degenerate end $E (= N\setminus M)$ and no parabolics. Let $S = E\cap M$ be the 
single boundary component of $M$.
Let $\EE = \{ E_\alpha \}$ denote the lifts of
$E$ to $\til{N}$, $S_\alpha = E_\alpha \cap \til{M}$, and $\SSS = \{ S_\alpha \}$ be the lifts of $S$ to $\til M$.
Suppose that  $\partial i: \partial \til{M} \rightarrow \partial \til{N}$ denotes the CT map given by Theorem \ref{red}.
Let $\partial{i} (a) = \partial{i} (b)$ for $a\neq b \in \partial \til{M}$. Let $\eta \subset \til{M} $ be the bi-infinite geodesic in $\til M$
joining $a, b$. Let $a_n \rightarrow a$ and $b_n \rightarrow b$ be points on $\eta$. Let $\eta_n$
(resp. $\rho_n$) be the geodesic in $\til M$ (resp. $\til N$) joining $a_n, b_n$.
By the converse direction of Lemma \ref{ea-strong}, we can approximate $\rho_n$ uniformly by  an electro-ambient quasigeodesic $\xi_n$
with respect to $(\til{N}, \EE)$. 

We now "reverse-engineer" an electro-ambient quasigeodesic $[a_n, b_n]_q$ with respect to 
 $(\til{M}, \EE)$ from $\xi_n$ as follows. This step is exactly the opposite of the corresponding step in the sketch before Theorem \ref{red}. 
 We replace any maximal segment of $\xi_n$ lying inside an $E_\alpha$
by a geodesic in the corresponding $S_\alpha \in \SSS$ to construct  $[a_n, b_n]_q$. Also, 
$[a_n, b_n]_q$ coincides with $\xi_n$ outside the $E_\alpha$'s.
By Lemma \ref{ea-strong}, $\eta_n$ lies in a uniformly bounded neighborhood of $[a_n, b_n]_q$. Also, since $[a_n, b_n]_q \setminus \cup \bigcup_\alpha
S_\alpha$ coincides with $\xi_n \setminus \cup \bigcup_\alpha E_\alpha$ and since $\xi_n$ converges to 
$\partial{i} (a) = \partial{i} (b)$ as $n \rightarrow \infty$, it follows that $[a_n, b_n]_q$ converges (in the Hausdorff metric on 
the Gromov compactification $\bbar{M}$) to a collection $\cup_r (c_r, d_r)$ of bi-infinite geodesics, with $\partial i (c_r)
= \partial i (d_r) =\partial{i} (a) = \partial{i} (b)$ for all $r$. By construction of $[a_n, b_n]_q$, each $(c_r, d_r)$
lies entirely in a single $S_\alpha$ and the CT map $\partial i_\alpha : S_\alpha \rightarrow E_\alpha$ identifies $c_r, d_r$.
This shows that the equivalence relation given by the CT map $\partial i: \til{M} \rightarrow \til{N}$ is generated by the 
equivalence relation given by the CT maps $\partial i_\alpha : S_\alpha \rightarrow E_\alpha$.
Corollary \ref{ptpre-relhyp}  below generalizes this argument to the relatively hyperbolic setup.

\begin{cor}\label{ptpre-relhyp} 
 Let $Y,X, \HH, \JJ, \BB = \{B_\alpha\}, \FF = \{F_\alpha\}, {\FF}_\alpha = \{F_{\alpha \beta}\}$ be as in Definition \ref{ei}
and $i:Y \rightarrow X$ be an {\bf ends-inclusion} of spaces. 

Also, let $\partial i: \partial \GG (Y, \HH ) \rightarrow \partial\GG (X, \JJ )$ be the induced Cannon-Thurston map on relative
hyperbolic boundaries as in theorem \ref{red}.
 Then $\partial i(a) = \partial i(b)$ for $a \neq b \in \partial \GG (Y, \HH )$ if and only if $a \RR b$ where
$\RR$ is the smallest closed equivalence relation generated by
the  CT relations $\RR_\alpha$ induced by $i_\alpha$. \end{cor}

\begin{proof} Let $\RR_Y$ denote the CT equivalence relation on $\partial Y$ induced by the CT map $\partial Y \rightarrow \partial X$
given by Theorem \ref{red}. We have to show that $\RR=\RR_Y$.

Since $i_\alpha: B_\alpha \rightarrow Y$ is a quasi-isometric embedding, it follows that $\RR_\alpha \subset \RR_Y$.
Hence the transitive closure $\RR_t$ of the union $\bigcup_\alpha \RR_\alpha$ is also contained in $\RR$. Finally, since
$\partial i:  \partial \GG (Y, \HH ) \rightarrow \partial\GG (X, \JJ )$ is continuous, it follows that $\RR_Y$ is a closed relation.
Hence $\RR \subset \RR_Y$.

It remains to show that $\RR_Y \subset \RR$.
Suppose  that $(a,b) \in \RR_Y$, i.e. $\partial i(a) = \partial i(b)$ for some $a \neq b \in \partial \GG (Y, \HH )$.
Then the geodesic $\eta= (a,b) \subset  \GG (Y, \HH )$ satisfies the following: \\
$\bullet$ If $a_n, b_n \in (a,b) \subset \GG (Y, \HH )$ are such that $a_n \rightarrow a, b_n \rightarrow b$, $\eta_n$ is the subsegment
of $\eta$ joining $a_n, b_n$, and $\rho_n$ is the geodesic in $\GG (X, \JJ )$
joining $a_n, b_n$, then $d_{\GG (X, \JJ )}( i(y), \rho_n) \rightarrow \infty$ as $n \rightarrow \infty$.

By the converse direction of Lemma \ref{ea-strong}, there exists $P \geq 1$ (independent of $a, b, n, a_n, b_n$) 
and  hyperbolic $P-$ quasigeodesic paths
$\xi_n$ such that 
\begin{enumerate} 
\item $\xi_n$ is an electro-ambient $P-$quasigeodesic with respect to $(X, \FF)$ lying within a hyperbolic distance $P$ of $\rho_n$.
\item There exists a sequence of points $a_n=a_{n,0}, a_{n,1}, \cdots , a_{n,k_n}=b_n$ on $\xi_n$ such that \\
a) The "odd subpaths" $\bbar{a_{n,2j},a_{n,2j+1}}$   of $\xi_n$ joining $a_{n,2j}$ and $a_{n,2j+1}$ have interiors disjoint from all $\GG(F_\alpha,
F_{\alpha \beta})$.\\
b) The "even subpaths" $\bbar{a_{n,2j+1},a_{n,2j+2}}$ of $\xi_n$ joining $a_{n,2j+1}$ and $a_{n,2j+2}$ are entirely contained in some $\GG(F_j,
F_{j \beta})$.
\end{enumerate}

By Remark \ref{geods-ei}, the odd subpaths of $\xi_n$ joining $a_{n,2j},a_{n,2j+1}$ are actually  $P-$quasigeodesics in $ \GG (Y, \HH )$.
Also, since $d_{\GG (X, \JJ )}( i(y), \rho_n) \rightarrow \infty$ as $n \rightarrow \infty$, it follows that 
$d_{\GG (Y, \HH )}( y, \bbar{a_{n,2j},a_{n,2j+1}}) \rightarrow \infty$ as $n \rightarrow \infty$.
In particular, for all $j$, $d_{\GG (Y, \HH )}( y,a_{n,j}) \rightarrow \infty$ as $n \rightarrow \infty$.

We shall now  reverse the construction used in the proof of Theorem \ref{red}.
Replace the even subpath $\bbar{a_{n,2j+1},a_{n,2j+2}}$ (of $\xi_n$) contained in $\GG(F_j,
F_{j \beta})$ by a geodesic $[a_{n,2j+1},a_{n,2j+2}]_q$ in the corresponding $\GG(B_j,
B_{j \beta})$. Interpolating the odd subpaths of $\xi_n$, we obtain an electro-ambient $P-$ quasigeodesic with respect to
 $((\GG(Y, \HH), \GG(B_\alpha,\{ B_{\alpha \beta} \}))$. Let $[a_n,b_n]_q$ denote this electro-ambient $P-$ quasigeodesic.

By Lemma \ref{ea-strong}, the geodesic in $\GG(Y, \HH)$ joining $a_n, b_n$ lies in a $K(=K(P))-$ neighborhood of $[a_n,b_n]_q$.
Passing to a subsequence if necessary, let $[a_n,b_n]_q$ converges to $[a,b]_q$ in the Hausdorff topology on 
the Gromov compactification $\bbar{\GG (Y, \HH )}$ of $\GG (Y, \HH )$. Let $(a,b)_q = [a,b]_q \cap \GG (Y, \HH )$. Then the geodesic $\eta$ lies in 
a $K-$ neighborhood of $(a,b)_q$. Thus  $(a,b)_q$ is a countable union of bi-infinite geodesics $(c_r,d_r) \subset \GG (Y, \HH )$, such that
$\eta$ lies in a $K-$ neighborhood of $\cup_r(c_r,d_r)$.  Here
$c_r, d_r \in \partial \GG (Y, \HH )$. Also, each such $(c_r,d_r)$ is a limit of geodesic segments contained in 
some (sequence of) $\GG(B_\alpha,\{ B_{\alpha \beta} \})$. Hence (passing to a further subsequence if necessary) we can assume that
each  $(c_r,d_r)$ is contained in some $\GG(B_r,\{ B_{r \beta} \})$.

Again, since each $c_r$ or $d_r$ is a limit of points ($a_{n,r}$) on $\xi_n$ and since 
$d_{\GG (X, \JJ )}( i(y), \xi_n)$ $ \rightarrow \infty$ as $n \rightarrow \infty$, 
all the $c_r, d_r$ get identified  with $a, b$  under the CT map $\partial i$.

Further, by the construction of $[a_n,b_n]_q$ from $\xi_n$, it follows that for any $r$, the pair 
$c_r, d_r$ get identified with each other under the CT map $\partial i_r$ (corresponding to the
inclusion of $\GG(B_r,\{ B_{r \beta} \})$ in $\GG(F_r,\{ F_{r \beta} \})$). Hence $(c_r, d_r) \in \RR_r$, where $\RR_r$ is the CT-relation induced
by $\partial i_r$.

Finally, since $\eta$ lies in a $K-$ neighborhood of $\cup_r(c_r,d_r)$, it follows that the pair $(a,b)$ is contained in the
smallest closed relation on $\partial  \GG (Y, \HH ) \times \partial  \GG (Y, \HH )$ generated by $\RR_r$, i.e. $(a,b) \in \RR$.
Hence $\RR_Y \subset \RR$ and the proof is complete.
\end{proof}

\subsection{Kleinian Groups with no Accidental Parabolics}

The first application of Theorem \ref{red} is to prove the existence of Cannon-Thurston maps for pared 3-manifolds with incompressible boundary and no accidental parabolics.
We  recall 
 the main Theorem of \cite{mahan-split}. Let $S$ be a complete finite area hyperbolic surface with fundamental group $H$. Nontrivial elements
of $H$ represented by peripheral loops of $S$ are called {\it parabolic elements} of $H$. Let $\til S$ denote the universal cover
of $S$. Note that $\til S$ is isometric to ${\mathbb{H}}^2$. Let $\bbar S = \til S \cup S^1 $ denote the Gromov compactification of
$\til S$. For $\rho$ a discrete faithful representation of $H$ into $Isom({\mathbb{H}}^3)$ taking parabolics to parabolics,
$\Gamma = \rho (H)$ is called a surface Kleinian group. If, in addition, $\rho$ does not send any non-parabolic element of $H$ to a 
parabolic, then $\Gamma$ is a surface Kleinian group without accidental parabolics.
In Theorem \ref{split} below, the convex core of ${\mathbb{H}}^3 / \Gamma$
will be denoted by $M$ and the union of $\til M$ with its limit set will be denoted by $\bbar M$. We are now ready to recall 
the main Theorem of \cite{mahan-split}:

\begin{theorem} \cite{mahan-split}
Let $\rho$ be a representation of a surface group $H$ (corresponding
to the surface $S$) into
$Isom({\mathbb{H}}^3)$ without accidental parabolics. Let $M$ denote (the convex
core of) ${\mathbb{H}}^3 / \rho
(H)$.  Further suppose that $i: S \rightarrow M$, taking
 parabolics to parabolics, induces a homotopy
equivalence.
  Then the inclusion
  $\tilde{i} : \widetilde{S} \rightarrow \widetilde{M}$ extends continuously
  to a map of the compactifications
  $\bbar{i} : \bbar{ S} \rightarrow \bbar{ M}$.  \label{split} \end{theorem}

\begin{theorem}
Suppose that $N^h \in H(M,P)$ is a hyperbolic structure 
on a pared manifold $(M,P)$ with no accidental parabolics. Further suppose that $N^h$ is not a doubly degenerate manifold. Let
$M_{gf}$ denotes a geometrically finite hyperbolic structure adapted
to $(M,P)$, then the map  $i: \widetilde{M_{gf}}
\rightarrow \widetilde{N^h}$ extends continuously to  a map of the compactifications
$\bbar{i}:  \bbar{{M_{gf}}}
\rightarrow \bbar{{N^h}}$.
\label{main}
\end{theorem}

\begin{proof} We first show that the lift of each end to $\til{N^h}$ is relatively quasiconvex. Suppose not.

Then by Lemma \ref{canary-qc}  a lift 
$\til E$ of an end of $N^h$ to $\til{N^h}$ is not relatively quasiconvex in $\til{N^h}$ if and only if there is a 
component $F$  of $\delta_0 M$ such that \\
\begin{enumerate}
 \item $F$ bounds a degenerate end other than $E$, and
\item $F$ is homotopic into $E$.
\end{enumerate}

If  $F$ is isotopic to a proper subsurface of $\partial E$, then
the boundary curves of $F$ necessarily have to be accidental parabolics contadicting the hypothesis.  

Else  $F$ is isotopic to all of $\partial E$, forcing $N^h$ is to be a doubly degenerate manifold
and again contadicting the hypothesis.  

Hence the map  $i: \widetilde{M_{gf}}
\rightarrow \widetilde{N^h}$ is an ends-inclusion.

The Theorem is now immediate consequence of Theorems \ref{red} and \ref{split}.
\end{proof}

\begin{rmk} \label{mainrmk} For the proof of Theorem \ref{main} to work it suffices to assume that each {\bf augmented} end 
of ${N^h}$ is relatively quasiconvex.
This will be useful in the next subsection when we deal with accidental parabolics.
\end{rmk}

To state the next theorem describing the point-preimages of the CT map, we set up some notation. Let $N$ be (the convex core of) a
hyperbolic structure 
on a pared manifold $(M,P)$ with relative Scott core $M_{gf}$. Let $\EE = \{ E_\alpha \}$ denote the  lifts of the (relative) ends of $N$
(i.e. the components of $N \setminus M_{gf}$). Let $\SSS_\alpha = E_\alpha \cap \til{M_{gf}}$.
  Let $\LL_\alpha$ denote the lift of the ending lamination (for the end corresponding to $E_\alpha$) to $\SSS_\alpha$. Each $\LL_\alpha$
induces an equivalence relation $\RR_\alpha$ on $\partial \til{M_{gf}}$ as follows: \\
$a \RR_\alpha b$ if and only if $a, b$ are ideal end-points of a leaf or complementary ideal polygon of $\LL_\alpha$. Let $\RR$ be the smallest 
closed equivalence relation (with respect to the product topology on $\partial \til{M_{gf}} \times \partial \til{M_{gf}}$) containing all the
equivalence relations $\RR_\alpha$.

In \cite{mahan-elct}  we also identify the point pre-images of the Cannon-Thurston map.
\begin{theorem} \cite{mahan-elct}
Let $G$ be a simply degenerate surface Kleinian group without accidental parabolics. Then the Cannon-Thurston map 
$\partial i: \partial \til{S} \rightarrow \partial \til{M}$ from the 
(relative) hyperbolic boundary of $G$ (which is the same as $\partial \til{S}$) 
to its limit set identifies precisely the end-points of leaves of the ending laminations.
More precisely, let $\RR$ denote the equivalence relation
on $\partial \til{S}$ given by $a\RR b$ iff $a, b$ are endpoints of a (lift of a) leaf of the ending
lamination or ideal boundary points of a complementary ideal polygon.
Then $\partial i(a) =  \partial i(b)$ if and only if $a \RR b$. \label{ptpre} \end{theorem}
Now, combining Theorem \ref{main},  Corollary \ref{ptpre-relhyp} and Theorem \ref{ptpre} we get:
\begin{theorem}
Suppose that $N^h \in H(M,P)$ is a hyperbolic structure 
on a pared manifold $(M,P)$ such that that  $N^h$ has no accidental parabolics.
 Let
$M_{gf}$ denotes a geometrically finite hyperbolic structure adapted
to $(M,P)$. Let  $\partial {i}:  \partial {\til{M_{gf}}}
\rightarrow \partial {\til{N^h}}$ be the Cannon-Thurston map extending  $i: \widetilde{M_{gf}}
\rightarrow \widetilde{N^h}$. 
Then $\partial i (a) = \partial i (b)$ for $a \neq b$ if and only if $(a,b) \in \RR$, where $\RR$ is the smallest 
closed equivalence relation containing the equivalence relations $\RR_\alpha$.
\label{ptprecor} \end{theorem}

\subsection{Accidental Parabolics} We shall now proceed to  remove the
restriction on accidental parabolics from Theorem \ref{main}. The proof proceeds by applying Theorems \ref{red} and \ref{split} twice successively.

\begin{theorem}
Suppose that $N^h \in H(M,P)$ is a hyperbolic structure 
on a pared manifold $(M,P)$ with incompressible boundary $\partial_0 M$. Suppose further that $N^h$ is not doubly degenerate. Let
$M_{gf}$ denotes a geometrically finite hyperbolic structure adapted
to $(M,P)$. Then the map  $i: \widetilde{M_{gf}}
\rightarrow \widetilde{N^h}$ extends continuously to the boundary
$\bbar{i}: \bbar{M_{gf}}
\rightarrow \bbar{N^h}$. 

 Let  $\partial {i}:  \partial {\til{M_{gf}}}
\rightarrow \partial {\til{N^h}}$ be the resulting Cannon-Thurston map extending  $i: \widetilde{M_{gf}}
\rightarrow \widetilde{N^h}$. 
Then $\partial i (a) = \partial i (b)$ for $a \neq b$ if and only if $(a,b) \in \RR$, where $\RR$ is the smallest 
closed equivalence relation containing the equivalence relations  generated by lifts of the ending laminations to 
$\til{M_{gf}}$.
\label{main2}
\end{theorem}

\begin{proof} First note that by Theorem \ref{main} and Remark \ref{mainrmk}, the Theorem follows when each augmented end is relatively quasiconvex.
Next, by Lemma \ref{canary-qc}, it follows that an  augmented end $E$ of $N^h$ is not relatively quasiconvex
if and only if there is a 
component $F$  of $\delta_0 M$ such that 
\begin{enumerate}
 \item $F$ bounds a degenerate end other than $E$, and
\item $F$ is homotopic into $E$.
\end{enumerate}

We construct another hyperbolic structure $W^h \in H(M,P)$ as follows: \\
For each augmented end $E$ of $N^h$ that is not relatively quasiconvex, let $F(E, i), i=1, \cdots k_E$ be the collection of 
components  of $\delta_0 M$ satisfying the 2 conditions above. Replace the  degenerate end 
having  $F(E, i), i=1, \cdots k_E$ as boundary by a geometrically finite end. We repeat this for every 
augmented end  that is not relatively quasiconvex. The resulting hyperbolic structure is denoted by $W^h$. We identify $W^h$ with its convex core
for convenience, i.e. we excise the geometrically finite (flaring) ends.

Each augmented end $E$ of $W^h$ is now relatively quasiconvex. By Theorem \ref{main} and Remark \ref{mainrmk}, 
the map  $j: \widetilde{M_{gf}}
\rightarrow \widetilde{W^h}$ extends continuously to the boundary
$\bbar{j}: \bbar{M_{gf}}
\rightarrow \bbar{W^h}$.

Each 
$F(E, i)$ is parallel to a subsurface of $\partial E$ and hence no other degenerate end can have boundary parallel to a subsurface of $F(E, i)$
unless $N^h$ is doubly degenerate (excluded by hypothesis).
It follows that the augmented ends bounded by $F(E, i)$ are relatively quasiconvex in $N^h$.
Hence the inclusion $j_2: \widetilde{W^h}
\rightarrow \widetilde{N^h}$ is an ends-inclusion and,  by Theorem \ref{main} and Remark \ref{mainrmk},  extends continuously to the boundary
$\bbar{j_2}: \bbar{W^h}
\rightarrow \bbar{N^h}$.

Since $i = j_2 \circ j$, it follows that the map  $i: \widetilde{M_{gf}}
\rightarrow \widetilde{N^h}$ extends continuously to the boundary
$\bbar{i}: \bbar{M_{gf}}
\rightarrow \bbar{N^h}$.

The last statement follows from (applying twice) the structure of the Cannon-Thurston map given by Theorem \ref{ptprecor}.
\end{proof}

\bibliography{ctredn}
\bibliographystyle{alpha}

\end{document}